\documentclass[11pt]{article}
\usepackage{amsfonts}
\usepackage{amssymb}
\title{\bf Kissing numbers -- a survey}
\author{
\begin{tabular}{cc} Peter Boyvalenkov, Stefan Dodunekov & Oleg Musin\thanks{This research is supported by the Russian government project 11.G34.31.0053, RFBR  grant 11-01-00735 and NSF grant DMS-1101688.} 
 \\
Inst. of Mathematics and Informatics & Department of Mathematics \\ 
Bulgarian Academy of Sciences & University of Texas at Brownswille \\ 
8 G. Bonchev str., 1113 Sofia, Bulgaria & 80 Fort Brown, TX 78520, USA
\end{tabular}} \textwidth 13.5cm \textheight 18.5cm \topmargin 0in
\def\abstractname{} 
\def\abstract{\if@twocolumn
\section*{\abstractname}
\else \small
\quotation  \noindent \fi}
\def\endabstract{\if@twocolumn\else\endquotation\fi}

\date{}

\begin{document}
\maketitle

\begin{abstract}\noindent{\sc Abstract.} The maximum possible
number of non-overlapping unit spheres that can touch a unit
sphere in $n$ dimensions is called kissing number. The problem for
finding kissing numbers is closely connected to the more general
problems of finding bounds for spherical codes and sphere
packings. We survey old and recent results on the kissing numbers
keeping the generality of spherical codes.
\end{abstract}

\section{Introduction}

How many equal billiard balls can touch (kiss) simultaneously
another billiard ball of the same size? This was the subject of a
famous dispute between Newton and Gregory in 1694. The more
general problem in $n$ dimensions, how many non-overlapping
spheres of radius 1 can simultaneously touch the unit sphere ${\bf
S}^{n-1}$, is called the kissing number problem. The answer
$\tau_n$ is called kissing number, also Newton number (in fact,
Newton was right, without proof indeed, with his answer
$\tau_3=12$) or contact number.

Further generalization of the problem leads to investigation of
spherical codes. A spherical code is a non-empty finite subset of
${\bf S}^{n-1}$. Important parameters of a spherical code $C
\subset {\bf S}^{n-1}$ are its cardinality $|C|$, the dimension
$n$ (it is convenient to assume that the vectors of $C$ span ${\bf
R}^n$) and the maximal inner product \[ s(C)=\max \{\langle x,y
\rangle : x,y \in C, x \neq y\}. \]

The function \[ A(n,s)=\max \{ |C|: \exists C \subset {\bf S}^{n-1} \mbox{
with } s(C) \leq s\} \] extends $\tau_n$ and it is easy to see that
$A(n,1/2)=\tau_n$. One also considers the function
\[ D(n,M)=\max \{ d(C)=\sqrt{2(1-s(C))}: \exists C \subset {\bf S}^{n-1} \mbox{ with }
|C|=M \} \] which is used in the information theory (cf.
\cite{CS,EZ1,Lev3}).

For $n \geq 3$ and $s>0$, only a few values of $A(n,s)$ are known.
In particular, only six kissing numbers are known: $\tau_1=2$,
$\tau_2=6$ (these two are trivial), $\tau_3=12$ (some incomplete
proofs appeared in 19th century and Sch\"utte and van der Waerden
\cite{SW} first gave a detailed proof in 1953, see also
\cite{Lee1,Was,Ans,Mus}), $\tau_4=24$ (finally proved in 2003 by
Musin \cite{Mus}), $\tau_8=240$ and $\tau_{24}=196560$ (found
independently in 1979 by Levenshtein \cite{Lev1} and
Odlyzko-Sloane \cite{OS}).

Note that Kabatiansky and Levenshtein have found an asymptotic upper bound $2^{0.401n(1+o(1))}$ for $\tau_n$ \cite{KL}. (Currently known the lower bound is $2^{0.2075n(1+o(1))}$ \cite{Wyn}.)

This survey deals with the above-mentioned values of $\tau_n$ and
mainly with upper and lower bounds in dimensions $n \leq 32$.
Some interesting advances during the last years are described.

Usually the lower bounds are obtained by constructions and the
upper bounds are due to the so-called linear programming
techniques and their extensions. We describe constructions which
often lead to the best known lower bounds. The upper bounds are
based on the so-called linear programming \cite{DGS,KL} and its
strengthening \cite{Mus,Pfe,Mus}. Applications were proposed by
Odlyzko-Sloane \cite{OS}, the first author \cite{Boy1}, and
strengthening by the third author \cite{Mus} and Pfender \cite{Pfe}.

Recently, the linear programming approach was strengthened as the
so-called semi-definite programming method was proposed by
Bachoc-Vallentin \cite{BV} with further applications by Mittelmann-Vallentin \cite{MV}.

\section{Upper bounds on kissing numbers}

\subsection{The Fejes T\'oth bound and Coxeter-B\"or\"oczky bound}

Fejes T\'oth \cite{FT1} proved a general upper bound on the
minimum distance of a spherical code of given dimension and
cardinality. In our notations, the Fejes T\'oth bound states that
\begin{equation}
\label{FTbound} D(n,M) \leq d_{FT} = \left(4-\frac{1}{\sin^2
\varphi_M} \right)^{1/2}
\end{equation}
where $\varphi_M=\frac{\pi M}{6(M-2)}$. This bound is attained for
$M=3,4,6$, and $12$. This gives four exact values of the function
$D(n,M)$ (but not necessarily implying exact values for $A(n,s)$).

First general upper bounds on the kissing numbers were proposed by
Coxeter \cite{Cox} and were based on a conjecture that was proved
later by B\"or\"oczky \cite{Bor}. Thus it is convenient to call
this bound the Coxeter-B\"or\"oczky bound.

Let the function $F_n(\alpha)$ be defined as follows:
\[ F_0(\alpha)=F_1(\alpha)=1, \]
\[ F_{n+1}(\alpha)=\frac{2}{\pi} \int_{(1/2)\arccos(1/n)}^{\alpha}
F_{n-1}(\beta(t)) dt \]
for $n \geq 1$, where
$\beta(t)=\frac{1}{2}\arccos \frac{\cos 2t}{1-2\cos 2t}$. This
function was introduced by Schl\"afli \cite{Sch} and is usually
referred to as Schl\"afli function.

In terms of the Schl\"afli function the Coxeter-B\"or\"oczky bound
is
\begin{equation}
\label{CBbound} A(n,s) \leq A_{CB}(n,s) =
\frac{2F_{n-1}(\alpha)}{F_n(\alpha)},
\end{equation}
where $\alpha = \frac{1}{2} \arccos \frac{s}{1+(n-2)s}$.

The bounds $\tau_n \leq A_{CB}(n,1/2)$ are weaker than the linear
programming bound to be discussed below. On the other hand, we
have
\[ A(4,\cos \pi/5)=120=A_{CB}(4,\cos \pi/5)=\frac{2F_3(\pi/5)}{F_4(\pi/5)} \]
(the lower bound is ensured by the 600-cell). The value $A(4,\cos
\pi/5)=120$ can be found by linear programming as well \cite{And}.
This suggests that the Coxeter-B\"or\"oczky bound can be better
than the linear programming bounds when $s$ is close to 1.

\subsection{Pure linear programming bounds}

The linear programming method for obtaining bounds for spherical
codes was built in analogy with its counter-part for codes over
finite fields which was developed by Delsarte \cite{Del}.
Delsarte-Goethals-Seidel \cite{DGS} proved in 1977 the main
theorem and it was generalized by Kabatianskii-Levenshtein
\cite{KL} in 1978.

The Gegenbauer polynomials \cite{AS,Sze} play important role in
the linear programming. For fixed dimension $n$, they can be
defined by the recurrence $P_0^{(n)}=1$, $P_1^{(n)}=t$ and
\[ (k+n-2)P_{k+1}^{(n)}(t)=(2k+n-2)tP_k^{(n)}(t)-kP_{k-1}^{(n)}(t)
                \mbox{ for } k \geq 1. \]
If \[ f(t)=\sum_{i=0}^m a_i t^i \] is a real polynomial, then
$f(t)$ can be uniquely expanded in terms of the Gegenbauer
polynomials as \[ f(t) = \sum_{k=0}^m f_kP_k^{(n)}(t). \] The
coefficients $f_i$, $i=0,1,\ldots,k$, are important in the linear
programming theorems.

\medskip

{\bf Theorem 1.} (Delsarte-Goethals-Seidel \cite{DGS},
Kabatianskii-Levenshtein \cite{KL}) {\it Let $f(t)$ be a real
polynomial such that

(A1) $f(t) \leq 0 $ for $-1 \leq t \leq s$,

(A2) The coefficients in the Gegenbauer expansion $f(t)=
\sum_{k=0}^m f_k P_k^{(n)}(t)$ satisfy $f_0>0$, $f_k \geq 0$ for
$i=1,\ldots,m$.

Then  $A(n,s) \leq f(1)/f_0$.}

\medskip

There are two cases, in dimensions eight and twenty-four, where
only technicalities remain after Theorem 1. The lower bounds
$\tau_8 \geq 240$ and $\tau_{24} \geq 196560$ are obtained by
classical configurations and the upper bounds are obtained by the
polynomials
\[ f_6^{(8,0.5)}(t)=(t+1)(t+\frac{1}{2})^2t^2(t-\frac{1}{2}) \]
and
\[ f_{10}^{(24,0.5)}(t)=(t+1)(t+\frac{1}{2})^2(t+\frac{1}{4})^2t^2
(t-\frac{1}{4})^2(t-\frac{1}{2}) \] respectively (the notations
will become clear later). Indeed, one may easily check that these two polynomial satisfy the
conditions (A1) and (A2) for the corresponding values of $n$ and
$s$ and therefore $\tau_8 \leq \frac{f_6^{(8,0.5)}(1)}{f_0}=240$
and $\tau_{24} \leq \frac{f_{10}^{(24,0.5)}(1)}{f_0}=196560$.

Together with the Gegenbauer polynomials we consider their
adjacent polynomials which are Jacobi polynomials $P_k^{(\alpha,\beta)}(t)$
with parameters
\[ (\alpha,\beta)=(a+\frac{n-3}{2},b+\frac{n-3}{2}) \]
where $a,b \in \{0,1\}$ (the Gegenbauer polynomials are obtained
for $a=b=0$). Denote by $t_k^{a,b}$ the greatest zero of the
polynomial $P_k^{(\alpha,\beta)}(t)$. Then \[ t_{k-1}^{1,1} <
t_k^{1,0} < t_k^{1,1} \] for every $k \geq 2$.

Denote
\[ \mathcal{I}_m =  \left\{
  \begin{array}{ll}
    \left[t_{k-1}^{1,1}, t_k^{1,0} \right], & \mbox{if } m=2k-1, \\[2pt]
    \left[t_k^{1,0},t_k^{1,1} \right],      & \mbox{if } m=2k,
  \end{array} \right.
\]
for $k=1,2,...$ and $\mathcal{I}_0 = [-1,t_1^{1,0})$.

Then the intervals $\mathcal{I}_m$ are consecutive and
non-overlapping. For every $s \in \mathcal{I}_m$, the polynomial
\[ f_m^{(n,s)}(t) = \left\{
    \begin{array}{ll}
       (t-s) \left( T_{k-1}^{1,0}(t,s) \right)^2,     & \mbox{if } m=2k-1, \\
       (t+1)(t-s)\left( T_{k-1}^{1,1}(t,s) \right)^2, & \mbox{if } m=2k,
    \end{array} \right.
\]
can be used in Theorem 1 for obtaining a linear programming bound.
Levenshtein \cite{Lev1} proved that the polynomials
$f_m^{(n,s)}(t)$ satisfy the conditions {\bf (A1)} and {\bf (A2)}
for all $s \in \mathcal{I}_m$. Moreover, all coefficients $f_i$,
$0 \leq i \leq m$, in the Gegenbauer expansion of $f_m^{(n,s)}(t)$
are strictly positive for $s \in \mathcal{I}_m$. Hence this
implies (after some calculations) the following universal bound.

\medskip

{\bf Theorem 2.} (Levenshtein bound for spherical codes
\cite{Lev1,Lev2}) {\it Let $n \geq 3$ and $s \in [-1,1)$. Then
\[ A(n,s) \leq  \left\{
\begin{array}{l}
   L_{2k-1}(n,s) = {k+n-3 \choose k-1}
                     \left[ \frac{2k+n-3}{n-1} -
                            \frac{P_{k-1}^{(n)}(s)-P_k^{(n)}(s)}{(1-s)P_k^{(n)}(s)}
                            \right] \\[2pt]
                    \hspace*{5cm} \mbox{for $s \in \mathcal{I}_{2k-1}$}, \\[4pt]
       L_{2k}(n,s) = {k+n-2 \choose k}
                  \left[ \frac{2k+n-1}{n-1} -
                         \frac{(1+s) \left( P_k^{(n)}(s)-P_{k+1}^{(n)}(s) \right)}
                              {(1-s) \left( P_k^{(n)}(s)+P_{k+1}^{(n)}(s)\right)} \right] \\[2pt]
                  \hspace*{5cm} \mbox{for $s \in \mathcal{I}_{2k}$.  }
\end{array} \right.
\]
}

In particular, one has $\tau_8 \leq L_6(8,1/2)=L_7(8,1/2)=240$ and
$\tau_{24} \leq L_{10}(24,1/2)=L_{11}(24,1/2)=196560$. The
Levenshtein bound can be attained in some other cases (cf. the
tables in \cite{Lev1,Lev2,Lev3}).

The possibilities for existence of codes attaining the bounds $L_m(n,s)$ were discussed in
\cite{BDL}. In particular, it was proved in \cite[Theorem 2.2]{BDL} that the even bounds $L_{2k}(n,s)$ can be only attained
when $s=t_k^{1,0}$ or $s=t_k^{1,1}$. This follows from a close investigation of the two-point distance distribution
\[ A_t=\frac{1}{|C|} \sum_{x \in C} |\{y \in C: \langle x,y \rangle=t\}|=\frac{1}{|C|}|\{ (x,y) \in C^2: \langle x,y \rangle=t\}| \]
of the possible $(n,L_{2k}(n,s),s)$ codes.

\medskip

On the other hand it was proved by Sidelnikov \cite{Sid} (see also
\cite[Theorem 5.39]{Lev3}) that the Levenshtein bounds are the
best possible pure linear programming bound provided the degree of
the improving polynomial is at most $m$. This restriction was
later extended by Boyvalenkov-Danev-Boumova \cite{BDD} to $m+2$
and the polynomials $f_m^{(n,s)}(t)$ are still the best.

However, in some cases the Leveshtein bounds are not the best
possible pure linear programming bounds. This was firstly
demonstrated in 1979 for the kissing numbers by Odlyzko-Sloane
\cite{OS}. Boyvalenkov-Danev-Boumova \cite{BDD} proved in 1996
necessary and sufficient conditions for existence of certain
improvements.

\medskip

{\bf Theorem 3.} \cite{BDD} {\it The bound $L_m(n,s)$ can be
improved by a polynomial from $A_{n,s}$ of degree at least $m+1$
if and only if $Q_j(n,s) < 0$ for some $j \geq m+1$. Moreover, if
$Q_j(n,s)<0$ for some $j \geq m+1$, then $L_m(n,s)$ can be
improved by a polynomial from $A_{n,s}$ of degree $j$.}

\medskip

For $s=1/2$ (the kissing number case) and $3 \leq n \leq 23$, $n
\neq 8$, the Levenshtein bounds are better that the
Coxeter-B\"or\"oczky  bounds but weaker than these which were
obtained by Odlyzko-Sloane \cite{OS}.

In three dimensions, the Levenshtein bound gives $\tau_3 \leq
L_5(3,1/2) \approx 13.285$ and it can be improved to $\tau_3 \leq
13.184$ which is, of course, not enough. Then Anstreicher
\cite{Ans} in 2002 and Musin \cite{Mus} in 2003 presented new
proofs which were based on strengthening the linear programming
and using spherical geometry on ${\bf S}^2$. The Musin's approach
will be discussed in more details below.

In four dimensions, we have $\tau_4 \leq L_5(4,1/2)=26$ and this
can be improved to $\tau_4 \leq 25.5584$ which implies that
$\tau_4$ is 24 or 25. Then Arestov-Babenko \cite{AB} proved in
2000 that the last bound is the best possible one can find by pure
linear programming. Earlier (in 1993), Hsiang \cite{Hsi} claimed a
proof that $\tau_4=24$ but that proof was not widely recognized as
complete. Musin \cite{Mus} presented his proof of $\tau_4=24$ in 2003
to finally convince the specialists. 

Odlyzko-Sloane \cite{OS} use discrete version of the condition (A1) 
and then apply the usual linear programming for $s=1/2$ and $3 \leq n \leq 24$. 
Their table can be seen in \cite[Chapter 1, Table 1.5]{CS}. 
Upper bounds for $25 \leq n \leq 32$ by linear programming were published in \cite{BD}. 
Now the first open case is in dimension five, where it is known that $40 \leq \tau_5 \leq 44$ 
(the story of the upper bounds is: $\tau_5 \leq L_5(5,1/2)=48$, $\tau_5 \leq 46.345$ from \cite{OS}),
$\tau_5 \leq 45$ from \cite{BV} and $\tau_5 \leq 44.998$ from \cite{MV}).

Let $n$ and $s$ be fixed, the Levenshtein bound gives $A(n,s) \leq
L_m(n,s)$ and it can be improved as seen by Theorem 3. In \cite{Boy3}, the first author 
proposed method for searching improving polynomials $f(t)=A^2(t)G(t)$, 
where $A(t)$ must have $m+1$ zeros in $[-1,s]$, $G(s)=0$ and $G(t)/(t-s)$ 
is a second or third degree polynomial which does not have zeros in $[-1,s]$. Moreover, 
one has $f_i=0$, $i=\in \{m,m+1,m+2,m+3\}$ for two or three consecutive coefficients in the Gegenbauer expansion
of $f(t)$. There restrictions leave several unknown parameters
which can be found by consideration of the partial derivatives of
$f(1)/f_0$ and numerical optimization methods. This approach was
realized (see \cite{Kaz}) by a programme SCOD. In fact, SCOD first checks for possible improvements by Theorem 3 and then
applies the above method. It works well for improving $L_m(n,s)$
for $3 \leq m \leq 16$ and wide range of $s$.

\subsection{Strengthening the linear programming}

The linear programming bounds are based on the following identity
\begin{equation}
  \label{main}
  |C|f(1)+\sum_{x,y\in C, x \neq y} f(\langle x,y\rangle)
      = |C|^2f_0 + \sum_{i=1}^k \frac{f_i}{r_i} \sum_{j=1}^{r_i}
        \left ( \sum_{x\in C} v_{ij}(x) \right )^2,
\end{equation}
where $C \subset {\bf S}^{n-1}$ is a spherical code,
\[ f(t)=\sum_{i=0}^k f_i P_i^{(n)}(t), \]
$\{ v_{ij}(x) : j=1,2,\ldots,r_i\}$ is an orthonormal basis of the
space $\mathrm{Harm}(i)$ of homogeneous harmonic polynomials of
degree $i$ and $r_i=\dim \,\mathrm{Harm}(i)$. In the classical
case (cf. \cite{DGS,KL} the sums of the both sides are neglected
for polynomials which satisfy (A1) and (A2) and this immediately
implies Theorem 1.

Musin \cite{Mus} strengthened the linear programming approach by
proposing the following extension of Theorem 1 which deals with a
careful consideration of the left hand side of (\ref{main}).

\medskip

{\bf Theorem 4.} \cite{Mus} {\it Let $f(t)$ be a real polynomial
such that

(B1) $f(t) \leq 0 $ for $t_0 \leq t \leq s$, where $-t_0>s$,

(B2) $f(t)$ is decreasing function in the interval $[-1,t_0]$,

(B3) The coefficients in the Gegenbauer expansion $f(t)=
\sum_{k=0}^m f_k P_k^{(n)}(t)$ satisfy $f_0>0$, $f_k \geq 0$ for
$i=1,\ldots,m$.

Then
\[ A(n,s) \leq \frac{\max \{ h_0,h_1,\ldots,h_{\mu}\}}{f_0}, \]
where $h_m$, $m=0,1,\ldots,\mu$, is the maximum of
$f(1)+\sum_{j=1}^m f(\langle e_1,y_j\rangle)$,
$e_1=(1,0,\ldots,0)$, over all configurations of $m$ unit vectors
$\{y_1,y_2,\ldots,y_m\}$ in the sphe\-rical cap (opposite of
$y_1$) defined by $-1 \leq \langle y_1,x \rangle \leq t_0$ such
that $\langle y_i,y_j\rangle \leq s$.}

\medskip

The proof of Theorem 4 follows from (\ref{main}) in a similar way
to the proof of Theorem 1 -- neglect the nonnegative sum in the
right hand side and replace the sum in the left hand side with its
upper bound \[ \sum_{i=1}^{|C|} \sum_{j: \langle y_i,y_j\rangle
\leq t_0} f(\langle y_i,y_j\rangle). \] Now observe that the last
expression does not exceed $\frac{\max \{
h_0,h_1,\ldots,h_{\mu}\}}{f_0}$.

Now the problems are to find $\mu$, choose $t_0$ and a polynomial
which mi\-ni\-mi\-zes the maximal value of
$h_0,h_1,\ldots,h_{\mu}$. In \cite{Mus} were found good polynomials
$f(t)$ by an algorithm which is similar to the algorithm of
Odlyzko-Sloane \cite{OS}. One easily sees that $h_0=f(1)$ and
$h_1=f(1)+f(-1)$. However, the calculation of the remaining
$h_m$'s usually requires estimations on
\[ S(n,M)=\min \max \{ s(C): C
\subset {\bf S}^{n-1} \mbox{ is a spherical code, } |C|=M\} \]
(cf. \cite{BDD,FT1,FT2,Lev3,Mus,SW}), observe that
$D(n,M)=\sqrt{2(1-S(n,M))}$). This approach was successfully
applied in dimensions three and four. In \cite{Mus} also noted
that this generalization does not give better upper bounds on the
kissing numbers in dimensions 5, 6, 7 and presumably can lead to
improvements in dimensions 9, 10, 16, 17, 18.

For $n=3$ and $s=1/2$ it is provedthat $\mu=4$, chooses
$t_0=-0.5907$ and finds suitable polynomial of degree 9 (similar
to these found in \cite{OS,Boy3,Kaz}) to show that $\tau_3=12$.
Analogously, for $n=4$ and $s=1/2$ he has $\mu=6$, $t_0=-0.608$
and certain polynomial of degree 9 to obtain $\tau_4=24$. The
calculations of $h_0$, $h_1$ and $h_2=\max_{\varphi \leq \pi/3}
\{f(1)+f(\cos \varphi)+f(-\cos(\pi/3-\varphi))\}$ are easy but
computations of $h_3,\ldots,h_6$ require numerical methods.

\subsection{Semidefinite programming}

Let $C= \{ x_i\} \subset \mathbb{S}^{n-1}$ be a spherical code, let $I \subset [-1,1)$ and let
\[ s_k(C,I) := \sum_{\langle x,y \rangle \in I, x,y \in C } \langle x,y \rangle^k =
|C| \sum_{t \in I} A_t t^k. \]
Odlyzko and Sloane \cite{OS} used in dimension 17 the constraints
\[ s_0(C,I_1) \leq |C|, \ s_0(C,I_2) \leq 2|C|, \]
where $I_1=[-1,-\sqrt{3}/2)$ and $I_2=[-1,-\sqrt{2/3})$, to improve on the LP bound. More general,
if it is know that the open spherical cap of angular radius $\varphi$ can contain at most $m$ points
of $C$ code with $S(C)=s$, where $\cos \varphi=t=\sqrt{s+(1-s)/(m+1)}$, then $s_0(C,I) \leq m|C|$, where $I=[-1,t)$.

Pfender \cite{Pfe} found the inequality
\[ s_2(C,I) \leq s_0(C,I)s+|C|(1-s), \]
where $I=[-1-\sqrt{s})$, and used it to improve the upper bounds for the kissing numbers in dimensions
9, 10, 16, 17, 25 and 26. In fact, the discussion in the preceding subsection can be 
viewed as in the following way: the third author \cite{Mus1,Mus2,Mus3,Mus} found a few inequalities for some linear
combinations of $s_k(C,I)$ for $0 \leq k \leq 9$, $s=1/2$ (the kissing numbers' case)
and certain $I=[-1,t_0]$, $t_0<-1/2$. In particular, that gave the proof that $\tau_4=24$ \cite{Mus}
and a new solution of the Thirteen spheres problem \cite{Mus2}.

This approach can be further generalized by consideration of the
three-point distance distribution
\[ A_{u,v,t}=\frac{1}{|C|} |\{ (x,y,z) \in C^3: \langle x,y \rangle=u, \langle x,z \rangle=v, \langle y,z \rangle=t\}| \]
(note that $A_{u,u,1}=A_u$). Here one needs to have $1+2uvt \geq u^2+v^2+t^2$. Bachoc-Vallentin  \cite{BV}
developed this to obtain substantial improvements for the kissing numbers in dimensions $n=4,5,6,7,9$ and 10.
Some numerical difficulties prevented Bachoc-Vallentin from furthers calculations but Mittelmann-Vallentin
\cite{MV} were able to overcome this and to report the best known upper bounds so far.

\section{Lower bounds on kissing numbers}

\subsection{Constructions A and B}

The idea for using error-correcting codes for constructions of
good spherical codes is natural for at least two reasons -- it
usually simplifies the description of codes and makes easier the
calculation of the code parameters. Leech-Sloane \cite{LS} make
systematic description of dense best sphere packings which can be
obtained by error-correcting codes and give, in particular, the
corresponding kissing numbers.

We describe Constructions A--B following \cite{CS}. Let $C$ be a
$(n,M,d)$ binary code. Then Construction A uses $C$ to build a
sphere packing in ${\bf R}^n$ by taking centers of spheres
$(x_1,x_2,\ldots,x_n)$, $x_i$ are integers, if and only if the
$n$-tuple \[ (x_1(\mbox{mod } 2),x_2(\mbox{mod }
2),\ldots,x_n(\mbox{mod } 2)) \] belongs to $C$.

The largest possible radius of nonoverlapping spheres is
$\frac{1}{2}\min\{2,\sqrt{d}\}$. The touching points on the sphere
with center $x$ are \[ 2^dA_d(x) \mbox{ if } d<4, \ 2n+16A_4(x)
\mbox{ if } d=4 \ 2n \mbox{ if } d>4, \] where $A_i(x)$ is the
numbers of codewords of $C$ at distance $i$ from $x$. Suitable
choices of codes for Construction A give good spherical codes for
the kissing number problem in low dimensions. The record lower
bounds for the kissing numbers which can be produced by
Construction A are shown in Table 1.

Let in addition all codewords of $C$ have even weight.
Construction B takes centers $(x_1,x_2,\ldots,x_n)$, $x_i$ are
integers, if and only if $(x_1(\mbox{mod } 2),x_2(\mbox{mod }
2),\linebreak \ldots,x_n(\mbox{mod } 2)) \in C$ and 4 divides
$\sum_{i=1}^n x_i$. The touching points on the sphere with center
$x$ are now
\[ 2^{d-1}A_d(x) \mbox{ if } d<8, \ 2n(n-1)+128A_8(x) \mbox{ if }
d=8, \ 2n(n-1) \mbox{ if } d>8. \] This, say complication, of
Construction A gives good codes for the kissing number problem in
dimensions below 24. It is remarkable that Construction B produces
the even part of the Leech lattice in dimension 24. The record
achievements of Construction B are also indicated in Table 1.

Having the sphere packings (by Constructions A and B, for example)
one can take cross-sections to obtain packings in lower dimensions
and can build up layers for packings in higher dimensions. This
approach is systematically used in \cite{CS} (see Chapters 5-8).

\subsection{Other constructions}

Dodunekov-Ericson-Zinoviev \cite{DEZ} proposed a construction
which develops the ideas from the above subsection by putting some
codes at suitable places (sets of positions in the original codes;
this is called concatenation in the coding theory). This
construction gives almost all record cardinalities for the kissing
numbers in dimensions below 24. Ericson-Zinoviev
\cite{EZ1,EZ2,EZ3} later proposed more precise constructions which
give records in dimensions 13 and 14 \cite{EZ3}.

\section{A table for dimensions $n \leq 32$}

The table of Odlyzko and Sloane \cite{CS,OS} covers dimensions $n
\leq 24$. Lower bounds by constructions via error-correcting codes
in many higher dimensions can be found in \cite{EZ1,CS} (see also
http://www.research.att.com/~njas/lattices/kiss.html). The
table below reflects our present (July 2012) knowledge in dimensions $n \leq 24$.

\begin{center}
\begin{tabular}{|c|c|c|}
  \hline
  Dimension  & Best known & Best known \\
     & lower bound & upper bound \\ \hline
  3 & 12 & 12 \\
  4 & 24 & 24 \\
  5 & 40 & 45 \\
  6 & 72 & 78 \\
  7 & 126 & 134 \\
  8 & 240 & 240 \\
  9 & 306 & 364 \\
  10 & 500 & 554 \\
  11 & 582 & 870 \\
  12 & 840 & 1357 \\
  13 & 1154 & 2069 \\
  14 & 1606 & 3183 \\
  15 & 2564 & 4866 \\
  16 & 4320 & 7355 \\
  17 & 5346 & 11072 \\
  18 & 7398 & 16572 \\
  19 & 10668 & 24812 \\
  20 & 17400 & 36764 \\
  21 & 27720 & 54584 \\
  22 & 49896 & 82340 \\
  23 & 93150 & 124416 \\
  24 & 196560 & 196560 \\
  \hline
\end{tabular}
\end{center}

The lower bounds in the Table above follow Table 1.5 from \cite{CS} apart from dimensions 13 and 14 taken from \cite{EZ3}.
The upper bounds are mainly taken from \cite{MV} (dimensions 5-7, 9-23).

Note that recently in \cite{CJKT} were found new kissing
configurations in 25 through 31 dimensions, which improve on the records set in 1982 by the laminated lattices.

\end{document}